# A Novel Decomposition Solution Approach for the Restoration Problem in Distribution Networks

Hossein Sekhavatmanesh, *Student Member, IEEE*, Rachid Cherkaoui, *Senior Member, IEEE*

*Abstract*—The distribution network restoration problem is by nature a mixed integer and non-linear optimization problem due to the switching decisions and Optimal Power Flow (OPF) constraints, respectively. The link between these two parts involves logical implications modelled through big-M coefficients. The presence of these coefficients makes the relaxation of the mixed-integer problem using branch-and-bound method very poor in terms of computation burden. Moreover, this link inhibits the use of classical Benders algorithm in decomposing the problem because the resulting cuts will still depend on the big-M coefficients. In this paper, a novel decomposition approach is proposed for the restoration problem named *Modified Combinatorial Benders (MCB)*. In this regard, the reconfiguration problem and the OPF problem are decomposed into master and sub problems, which are solved through successive iterations. In the case of a large outage area, the numerical results show that the MCB provides, within a short time (after a few iterations), a restoration solution with a quality that is close to the proven optimality when it can be exhibited.

## NOMENCLATURE

### A. Parameters

| | |
|---|---|
| $w_{re}, w_{sw}, w_{op}$ | Weighting factors of the objective function terms (p.u.) |
| $D_i$ | Importance factor of the load at bus $i$ (p.u.) |
| $\lambda_{ij}(\lambda_i)$ | The operation time of line switch $ij$/load breaker $i$ (hour). |
| $P_{i,t}^D$ ($Q_{i,t}^D$) | Active (/Reactive) power demand at bus $i$, at time $t$ (p.u.). |
| $r_{ij}$ ($x_{ij}$) | Resistance (Reactance) of line $ij$ (p.u.). |
| $v^{max}$ ($v^{min}$) | Maximum (Minimum) limits of voltage magnitude (p.u) |
| $f_{ij}^{max}$ | Maximum current flow rating of line $ij$ (p.u) |
| $P_{i,max}^{inj}$ | Active power capacity of resource at node $i$ (p.u.) |
| $S_{i,max}^{inj}$ | Apparent power capacity of resource at node $i$ (p.u.) |
| $M$ | A large multiplier |

### B. Variables

| | |
|---|---|
| $\mu_{ij}$ | Binary decision variable indicating if the line $ij$ equipped with a switch is energized or not (1/0) |
| $\alpha_{i,t}$ | Binary decision variable indicating if at time $t$ the load at node $i$ is supplied or rejected (1/0) |
| $\beta_{ij}$ | Continuous variable indicating if the line $ij$ is oriented from node $i$ to node $j$ or not. |
| $\psi_{ij}$ | Auxiliary flow that is travelling in line $ij$ from node $i$ to node $j$. |
| $\phi_i$ | Auxiliary variable indicating if the node $i$ is energized or not (1/0). |
| $F_{ij,t}$ | Square of current flow magnitude in line $ij$, at time $t$ (p.u) |
| $P_{ij,t}$ ($q_{ij,t}$) | Active (Reactive) power flow in line $ij$, starting from node $i$, at time $t$ (p.u.). |
| $P_{i,t}^{inj}$ ($Q_{i,t}^{inj}$) | Active (Reactive) power injection from the substation or DGs at node $i$, at time $t$ (p.u.). |
| $U_{i,t}$ | Square of voltage magnitude at bus $i$, at time $t$ (p.u.). |

### C. Indices

| | |
|---|---|
| $i,j$ | Index of nodes |
| $ij$ | Index of branches |
| $t$ | Index of time |
| $q$ | Index of iteration in the decomposition approach |

### D. Sets

| | |
|---|---|
| $N$ | Set of nodes |
| $N^*$ | Set of nodes in the off-outage area |
| $W$ | Set of lines (plus tie-lines) |
| $W^*$ | Set of lines (plus tie-lines) in the off-outage area |
| $W^S = W_{ava}^S \cup W_{int}^S \cup W_{sec}^S$ | Set of lines (plus tie-lines) in the off-outage area equipped with switches |
| $W_{ava}^S$ | Set of lines hosting available tie-switches |
| $W_{int}^S$ | Set of lines hosting internal tie-switches |
| $W_{sec}^S$ | Set of lines hosting sectionalizing switches |
| $\Omega_{DG}$ | Set of nodes hosting DGs |

## I. INTRODUCTION

After a failure in a radial distribution network, once the fault is isolated, the area downstream to the fault place remains unsupplied. This area is called *off-outage area*. The aim of the restoration operation is to restore the maximum energy of loads within this off-outage area while minimizing the total switching operation time [1]. In order to achieve this goal in passive distribution networks, the only possible action is to transfer the unsupplied loads to the healthy neighboring feeders (Fig. 1.). This reconfiguration is deployed through changing the status of normally-closed (sectionalizing) and normally-open (tie) switches. The tie-switches that are between the faulted feeder and a healthy feeder are called *available tie-switches* (T2, T3, and T4 in Fig. 1.). The tie-switches with both ending nodes inside the off-outage area are referred to as *internal tie-switches* (T5 in Fig. 1.). The resulting new configuration of the network remains for a so-called *restorative period* that starts from the fault isolation instant until the time when the faulted element is repaired. After the restorative period, the original configuration of the network will be restored.

The highly increasing penetration of Distributed Generators (DGs) in Active Distribution Networks (ADNs) introduces, among others, new restoration actions besides the switching operations. Among these additional restoration actions are DG power set point modifications. The incorporation of these actions could lead to a more efficient restoration solution.

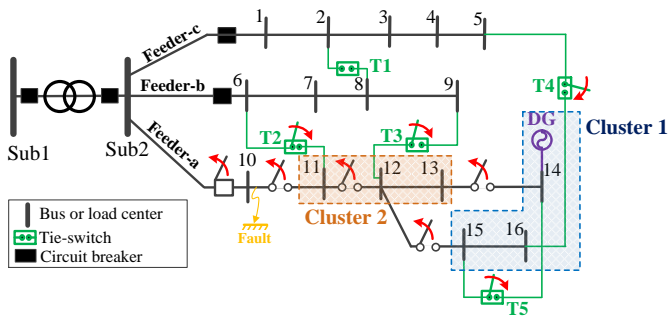

Fig. 1. A simple distribution network under post fault conditions.

However, this leads to increases the complexity of the restoration problem. Therefore, the challenge is to find an automatic and efficient solution for the restoration problem as an operator decision support in distribution networks while considering their new active status.

From mathematical formulation point of view, the restoration problem is an NP-hard combinatorial optimization problem. It is in the form of a mixed-integer and non-linear problem, respectively, due to the switching decisions and Optimal Power Flow (OPF) constraints. The OPF is a known challenging optimization problem. It has been the main building block for the formulation of many control, operation, and planning problems in power systems. In this respect, an AC-OPF should be integrated into the restoration problem in order to check the feasibility of the obtained solution concerning the technical constraints such as voltage and current limits. In order to deal with the non-polynomial hardness of such an OPF-based optimization problem, some researchers applied meta-heuristic methods as the solution approach. Among these methods are Genetic Algorithm (GA) [2], Particle Swarm Optimization (PSO) [3], Ant Colony [4], and Tabu Search [5]. The main problem with these methods is that they are in general computationally exhaustive and could fail to find an existing feasible and good-enough solution in a reasonable time complying with online restoration requirements.

In order to address this weakness, some researchers proposed to use heuristic approaches. These approaches use graph search methods to find a suitable topology of the network with regard to the restoration criteria [6], [7]. It means that the heuristic restoration algorithms will be applicable only for specific network topologies and cannot be generalized for any network topology. Moreover, as it is shown in [8], the restoration solution provided by these heuristic algorithms could be very far from the global optimal solution.

The mathematical analytical programming started to be an option for solving the restoration problem shortly after that some convex relaxation methods were proposed for the OPF problem. The authors of [9]-[10] used semidefinite programming (SDP) relaxation in order to solve it for radial unbalanced grids. Regarding the grid unbalancing, a distributed optimization technique is developed in [11] based on the Alternating Direction Method of Multipliers (ADMM). In this paper, we assume that the distribution network is operated in a radial and balanced fashion. In this respect, relaxation methods are proposed in the literature for the OPF problem in the form of either Mixed-Integer Quadratic Constraint Programming (MIQCP) [12], [13] or Mixed-Integer Second-Order Cone Programming (MISOCP) model [14]–[16]. The discussion given in Appendix I gives an insight about possible extension of the proposed restoration approach to unbalanced grids.

The convexity of the power flow formulation allows finding a solution with proven optimality up to the desired accuracy. However, the computation burden of the resulting optimization problem could be intractable depending on the dimension of the distribution network. This drawback inhibits the use of mathematical programming methods to tackle the restoration problem as a multi-period and/or multi-scenario optimization problem.

In order to address this problem, different studies proposed decomposition methods. Some of them account for the uncertainty of load demands and DG injections using stochastic or robust optimization [17]–[21]. They proposed decomposition algorithms where the innermost problem is assigned to find the worst realization of uncertain parameters given a fixed radial configuration. In the uppermost level, the deterministic restoration problem is solved while fixing the uncertain parameters to their worst-case realizations found in the inner level problem. For solving this decomposed restoration problem, different solution approaches have been applied. Among the most important ones are stochastic rolling-horizon optimization method [18], Information Gap Decision Theory (IGDT) [19], and column-and-constraint generation algorithm [20], [21].

The restoration strategy should consider the time-varying loads in order to provide a unique configuration valid throughout the whole restorative period. Actually, none of the above-mentioned papers accounts for the time-varying loads. If those decomposition algorithms are used for solving a such multi-period optimization problem, the deterministic restoration problem in the uppermost level will be computationally intractable. In this regard, papers [22], [23] propose to partition the time window of the optimization problem in clusters with similar load levels. Then, the reconfiguration problem is solved sequentially for each cluster of time instants. The weakness of this strategy is that the solutions at a given sequence are not influenced by the future sequences during the rest of the optimization process. Therefore, this approach cannot be applied to solve multi-period restoration problems (such as the one studied in this paper), where the feasible and optimal solution at one time step depends on the solution of the problem in the previous and next steps.

In order to relax the computation burden of the multi-period restoration problem, the authors in [17], [24], [25] propose to use a linear load flow formulation instead of the original and non-linear formulation. Since in the restoration strategy, the reconfigured network could be operated very close to its capacity envelop, applying the linear load flow model may result in an infeasible solution regarding the network safety constraints (e.g. voltage and current limits).

Another approach for relaxing the computation burden of the restoration problem is provided in [26]. According to this methodology, the restoration problem is solved in two stages. In the first stage, the post-restoration topology is determined using a heuristic approach. The set of loads to be restored and the outputs of sources are determined in the second stage while fixing the network topology to the one obtained in the first stage. The optimization problems in the first and second stages, which are referred as reconfiguration problem and load pickup problem, respectively, are mutually interdependent. The decoupling of these two problems as proposed in [26] could lead either to no feasible solution or to a solution very far from the optimal one.

In order to address all the afore-mentioned weaknesses, we propose a modification to the combinatorial Benders method so that it can be used for the multi-period restoration problem while considering case studies of realistic size. Combinatorial Benders method was firstly proposed by Hooker for solving optimization problems which include conditional constraints [27, p.]. This method has been used in many different applications such as in circuit verification problems [28], map labeling problem [29] and asymmetric travelling salesman problem [30].

In this paper, a two-stage analytical formulation is proposed for the restoration problem. It is solved thanks to a so called Modified Combinatorial Benders algorithm. Compared with the state-of-the-art, the major contributions of this paper are the following:

1- A novel decomposition approach is proposed for the restoration problem while considering inter-temporal constraints (e.g. varying loads) and control actions (e.g. DG power set points). The original problem is in the form of a multi-period, mixed integer, and non-linear optimization problem. Thanks to the proposed decomposition approach, the restoration problem is made tractable for analytical solvers in case of a grid of realistic size in a multi-period optimization problem.
2- The standard combinatorial Benders method is augmented with new cuts identifying binary variable combinations that are either infeasible or non-optimal. In this regard, the proposed cuts distil the search space of the optimization problem at a given iteration into a smaller subset that includes the global optimal solution. Therefore, compared with the standard combinatorial Benders, the proposed MCB approach converges in less number of iterations. The quality of this final solution is close to the global optimal solution.
3- A convex AC-power flow formulation is integrated into the decomposed formulation proposed for the restoration problem in order to accurately model the electrical operational constraints (e.g. voltage and current limits).

The remainder of this paper is organized as follows. The integrated mathematical model of the restoration problem is presented in section II. The decomposition approach is presented in section III. In this section, first, the details of the Combinatorial Bender decomposition method are explained. Then, the proposed modified one is provided. Section IV illustrates different case studies verifying the main advantages of the proposed decomposition algorithm with respect to the integrated analytical optimization method. Finally, section V concludes the paper highlighting the main contributions.

## II. PROBLEM FORMULATION

In this section, an integrated model of the time-dependent restoration formulation is presented in the form of a MISOCP problem. This problem encompasses three groups of decision variables namely, I) the binary variables $\mu_{ij}$ represent the energization status of line switch $ij$, II) the binary variables $\alpha_{i,t}$ account for the status of the load at node $i$ at time $t$, III) and the continuous variables $P_{i,t}^{DG}/Q_{i,t}^{DG}$ are associated to the active/reactive power set points of DG at node $i$ at time $t$. The targeted restoration strategy is a multi-period optimization problem in the sense that the decision for the load pickup and DG power set points varies with time. The line switching variables ($\mu_{ij}$) do not vary with time, since it is aimed to provide a single new network configuration for the whole restorative period.

Minimize: $F^{obj} = W_{re}.F^{re} + W_{sw}.F^{sw} + W_{op}.F^{op}$ (1.a)

where,

$$F^{re} = \sum_t \sum_{i \in N} D_i.(1 - \alpha_{i,t}).P_{i,t}^D$$

$$F^{sw} = \sum_{(i,j) \in W_{tie}^S} \mu_{ij}.\lambda_{ij} + \sum_{(i,j) \in W_{sec}^S} (1 - \mu_{ij}).\lambda_{ij}$$

$$F^{op} = \sum_t \sum_{(i,j) \in W} r_{ij}.F_{ij,t}$$

Subject to:

$$\begin{cases} 0 \leq \beta_{ij} \leq 1 & \forall(i,j) \in W^S \\ \beta_{ij} + \beta_{ji} = \mu_{ij}, & \forall(i,j) \in W^S \\ \beta_{ij} = \mu_{ij}, \beta_{ji} = 0, & \forall(i,j) \in W_{ava}^S \end{cases}$$ (1.b)

$$\begin{cases} \phi_i = \sum_{j:(i,j) \in W^*} \beta_{ji} \leq 1, & \forall i \in N^* \\ \phi_i = 1, & \forall i \in N \setminus N^* \end{cases}$$ (1.c)

$$\begin{cases} 0 \leq \Psi_{ij} \leq M.\beta_{ij} & \forall(i,j) \in W^S \\ 0 \leq \Psi_{ji} \leq M.\beta_{ji} & \forall(i,j) \in W^S \\ \sum_{\forall j^*:(j^*,i) \in W^S} (\Psi_{j^*i}) = \sum_{\forall j^*:(i,j^*) \in W^S} (\Psi_{ij^*}) + \phi_i & \forall i \in N^* \\ \sum_{\forall(i,j) \in W_{ava}^S} \Psi_{ij} = \sum_{\forall i \in N^*} \phi_i \end{cases}$$ (1.d)

$$\begin{cases} \alpha_{i,t} \leq \phi_i, & \forall i \in N^*, \forall t \\ \alpha_{i,t} = 1, & \forall i \in N \setminus N^*, \forall t \end{cases}$$ (1.e)

$$0 \leq \alpha_{i,t-1} \leq \alpha_{i,t} \leq 1 \qquad \forall i \in N \setminus N^*, \forall t \quad (1.f)$$

$$-M.(1 - \mu_{ij}) \leq U_{i,t} - U_{j,t} - 2(r_{ij}.p_{ij,t} + x_{ij}.q_{ij,t})$$
$$\leq M.(1 - \mu_{ij}) \qquad \forall ij \in W, \forall t \quad (1.g)$$

$$\begin{cases} p_{ij,t} = (\sum_{\substack{i^* \neq i \\ (i^*,j) \in W}} p_{ji^*,t}) + r_{ij}.F_{ij,t} + \alpha_{i,t}.P_{j,t}^D - P_{j,t}^{inj} \\ q_{ij,t} = (\sum_{\substack{i^* \neq i \\ (i^*,j) \in W}} q_{ji^*,t}) + x_{ij}.F_{ij,t} + \alpha_{i,t}.Q_{j,t}^D - Q_{j,t}^{inj} \end{cases}$$
$$\forall ij \in W, \forall t \quad (1.h)$$

$$\left\| \begin{matrix} 2p_{ij,t} \\ 2q_{ij,t} \\ F_{ij,t} - U_{i,t} \end{matrix} \right\|_2 \leq F_{ij,t} + U_{i,t} \qquad \forall(i,j) \in W, \forall t \quad (1.i)$$

$$\begin{cases} 0 \leq P_{i,t}^{inj} \leq P_{i,max}^{inj} \\ \left\| \begin{matrix} P_{i,t}^{inj} \\ Q_{i,t}^{inj} \end{matrix} \right\|_2 \leq S_{i,max}^{inj} \end{cases} \quad \forall i \in \Omega_{DG} \cap N, \forall t \quad (1.j)$$

$$\begin{cases} 0 \leq F_{ij,t} \leq \mu_{ij} \cdot f_{ij}^{max\,2} \\ -M.\mu_{ij} \leq p_{ij,t} \leq M.\mu_{ij} \\ -M.\mu_{ij} \leq q_{ij,t} \leq M.\mu_{ij} \end{cases} \quad \forall ij \in W, \forall t \quad (1.k)$$

$$v^{min\,2}.\phi_i \leq U_{i,t} \leq v^{max\,2}.\phi_i \quad \forall i \in N, \forall t \quad (1.l)$$

The objective function (1.a) tends to minimize the weighted total costs associated with the reliability ($F^{re}$), switching ($F^{sw}$), and operational ($F^{op}$) objectives, in decreasing order of priority. This hierarchical priority is enabled using $\epsilon$−constraint method [31]. The reliability cost is the total energy not supplied of the loads while accounting for their importance factors. The switching cost is formulated as the summation of two sub terms associated with the total operation time of tie-switches and sectionalizing switches, respectively. Finally, the least priority objective term is the operational term formulated as the total active power loss. As the formulation of the restoration problem is not the focus of this paper, only a very brief description of constraints (1) will be provided in the following. The readers are recommended to go through [8] for detailed explanation.

Constraints (1.b)-(1.d) model the reconfiguration problem ensuring the radial topology of the network [32]. Constraints (1.d) ensure that all the nodes are connected to the substation via flows of a virtual commodity represented by auxiliary and continuous variables $\Psi_{ij}$ and $\Psi_{ji}$. A brief explanation of these constraints is provided in Appendix IV. Using these constraints, we avoid to create an isolated loop in the off-outage area with the loads that are supplied in an islanded way using an existing DG in that loop. Constraints (1.e)-(1.f) formulate the load pickup problem. For an energized node $i$ in the off-outage area ($\phi_i = 1$), a decision is made in (1.d) with binary variable $\alpha_{i,t}$, indicating if its load is restored at time $t$ or rejected (1/0). As formulated in (1.f), it is assumed that once a given load is restored at a given time, no further interruption is permitted during the subsequent time slots of the restorative period. The relaxed formulation of AC-OPF is presented in (1.g)-(1.j). The aim of this part is to dispatch the active/reactive power set points of DGs, while respecting all the security constraints in the reconfigured network. The reconfiguration problem is linked to the AC-OPF problems in (1.k) and (1.l) using variables $\mu_{ij}$ and $\phi_i$, respectively. These links inhibit the use of classical Benders algorithm in decomposing the problem because the resulting cuts will still depend on the big-M coefficients used in (1.k) and (1.l).

In the subsequent discussion, it is aimed to present a tractable approach for solving the restoration problem in a multi-period optimization environment. In this respect, the following compact form of the restoration model will be used to represent the above extensive formulation.

$$P := \min_{x,y,z,u} \sum_{i \in N^*} f_i(1-x_i) + \sum_{ij \in W^*} e_{ij} y_{ij} + \sum_{ij \in W} d_{ij} z_{ij} \quad (2.a)$$

Subject to:
$$\mathbb{C}_1(x,y) \neq \emptyset \quad (2.b)$$
$$\mathbb{C}_2(x,y,z,u) \neq \emptyset \quad (2.c)$$

where:

$$\mathbb{C}_1(x,y) := \begin{cases} y \in \Gamma & (3.a) \\ Ax \geq a & (3.b) \\ By + Cx \geq b & (3.c) \end{cases}$$

$$\mathbb{C}_2(x,y,z,u) :=$$
$$\begin{cases} if \quad \eta_\sigma(y) = 1 \quad then \quad D_\sigma z \geq c_\sigma, \quad \forall \sigma = 1,2,\ldots,\sigma_{max} & (4.a) \\ Ez = Jx - Fu & (4.b) \\ \|G_l z\| \leq g_l^T z, & \forall l \in W & (4.c) \\ \|H_i u\| \leq h_i, & \forall i \in \Omega_{DG} & (4.d) \end{cases}$$

where, $y$ and $x$ are vectors of binary variables indicating, respectively, line switching variables ($\mu_{ij}$) and load pickup variables ($\alpha_{i,t}$). Continuous variables are represented by vectors $u$ and $z$, which refer, respectively, to the DG power set point variables ($P_{i,t}^{DG}/Q_{i,t}^{DG}$), and the rest of state variables related to the optimal power flow constraints at each time $t$ (such as $U_{i,t}, F_{ij,t}, \ldots$).

The three terms of (2.a) represent, respectively, the reliability, switching, and operational objective terms formulated in (1.a). $\mathbb{C}_1$ is expressed in (3) as the set of constraints only on the binary variables ($x$ and y). In (3.a), $\Gamma$ is the set of radial network configurations described by (1.b)-(1.d). Constraint (3.b) accounts for (1.f) as the load pickup formulation. Constraint (3.c) represents (1.e) as the link between reconfiguration and load pickup problems.

As given in (4), $\mathbb{C}_2$ represents the set of AC-OPF constraints which are linked to the binary variables (x and y). The link between the reconfiguration and AC-OPF problem is given in (4.a) in the form of conditional constraints, which are linearized using big-M formulation as formulated in (1.k) and (1.l). They mean that if a certain condition on $y$ variables holds ($\eta_\sigma(y) = 1$), a constraint on $z$ variables is added to the optimization problem[1]. Equation (4.b) accounts for the voltage equation and the power balance equation given in (1.g) and (1.h), respectively. Equations (4.c) and (4.d) represent second-order constraints associated with the current flow equation (1.i), and DG apparent power capacity limit (1.j). In (4.d), $h_i$ refers to the apparent power capacity of the DG at node $i$. The other notations used in (2),(3) and (4) are matrices (the ones in capital letter) or vectors (the ones in small letter) of parameters.

III. SOLUTION STRATEGY

In this section, a novel decomposition approach is proposed for the restoration problem named Modified Combinatorial Benders. In this regard, the reconfiguration problem and the OPF problem are decomposed into master and sub problems, which are solved through successive iterations.

---

[1] For example, if line $ij$ is energized ($\mu_{ij} = 1$), then the current flow in this line must be less than its ampacity limit ($F_{ij,t} \leq f_{ij}^{max\,2}$).

## A. Master Problem

In the outer level of the proposed decomposition strategy, the master problem is solved. The master problem is the same as the original problem $P$ while removing the AC-OPF constraints ($\mathbb{C}_2$). However, the operational constraints (e.g. voltage and current limits) are not completely disregarded as they are represented by DistFlow constraints ($\overline{\mathbb{C}_2}$).

$$\mathcal{M}^{(q)} := \min_{x,y,z,u} \sum_{i \in N^*} f_i(1-x_i) + \sum_{ij \in W^*} e_{ij} y_{ij} \quad (5.a)$$

Subject to:
$$\mathbb{C}_1(x,y) \neq \emptyset \quad (5.b)$$
$$\overline{\mathbb{C}_2}(x,y,z,u) \neq \emptyset \quad (5.c)$$

$$\begin{cases} \sum_{i \in \mathbb{X}(v,y^{(k)})} f_i(1-x_i) \geq \mathcal{L}_v(y^{(k)}) \; : y \in \mathbb{Y}(v,y^{(k)}), \mathcal{L}_v(y^{(k)}) \neq \emptyset \\ y \notin \mathbb{Y}(v,y^{(k)}) \qquad\qquad\qquad : \mathcal{L}_v(y^{(k)}) = \emptyset \end{cases}$$
$$\forall v = 1,2,\dots,v_{max}(k), k = 1,\dots,q-1 \quad (5.d)$$

where:
$$\overline{\mathbb{C}_2}(x,y,z,u) :=$$
$$\begin{cases} if \; \eta_\sigma(y) = 1 \; then \; D_\sigma z \geq c_\sigma, \quad \forall \sigma = 1,2,\dots,\sigma_{max} & (6.a) \\ Ez = Jx - Fu & (6.b) \\ \overline{G}_l^T z \leq \overline{s}_l, & \forall l \in \Omega_{Sub} & (6.c) \\ \overline{H}_i^T u \leq \overline{h}_i, & \forall i \in \Omega_{DG} & (6.d) \end{cases}$$

$\overline{\mathbb{C}_2}$ is formulated in (6) as the set of linear DistFlow constraints linked to the binary variables. Constraints (6.c) and (6.d) are the linearized formulation of line ampacity and DG apparent power capacity limits, where $\overline{s}_l$ and $\overline{h}_i$ induce relaxed line ampacity limit and DG apparent power limit, respectively. This linearization technique is according to the technique presented in [20]. The detailed formulation of the DistFlow constraints is given in [33].

Two sets of constraints (5.d) are denominated as optimality cuts and feasibility cuts, respectively. These constraints represent the modified version of the Combinatorial Benders cuts introduced in [34]. The main idea behind this modification is to have recourse functions providing information not only on the feasibility of each solution but also on its optimality. According to the standard Combinatorial Benders method, all the binary variables must be set as the complicating variables. They are determined by the master problem and the sub problem just to evaluate the feasibility of the solution. In case that the solution is infeasible, a cut will be added to the master problem for the next iteration to remove the corresponding set of infeasible binary variables from the solution space. In the proposed modified version, only a subset of binary variables is fixed in the sub problem and defined as the complicating variables. The other binary variables that are not fixed in the sub problem are called floating variables. In the formulation provided in (5) $y$ and x represent, respectively, complicating and floating variables. In the case where the solution of the master problem at iteration $k$ ($y^{(k)}$) leads to no feasible solution in the sub problem, the feasibility cuts are augmented by the second constraint of (5.d). If the sub problem at iteration k is feasible, its optimal solution $\mathcal{L}_v(y^{(k)})$ is used to augment the optimality cuts according to the first constraint of (5.d). The constraints given in (5.d) are non-convex and need to be linearized. This linearization together with the formulation of $\mathbb{Y}(v,y^{(k)})$ is provided in section III.C.

$$LB^{(q)} = \mathcal{M}^{(q)} \quad (7)$$

As given in (7), the master problem $\mathcal{M}^{(q)}$ formulated in (5.c) provides a lower bound for the original optimization problem $P$ expressed in (2). Actually, unlike in the case of AC-OPF formulation in the sub problem, we choose relaxed limits for the electrical operation constraints (e.g. voltage limits) of the DistFlow formulation in the master problem. Therefore, the feasible region of $\mathcal{M}^{(q)}$ under DistFlow constraints is relaxed in comparison to the feasible region of $P$ under AC-OPF constraints.

## B. Sub Problem

Once the optimal configuration is found in the master problem, the next step is to find the optimal load pickup solution, if any, for the obtained configuration subject to the AC-OPF constraints. When we fix the network configuration, the topology of the off-outage area will be partitioned accordingly into several clusters. A cluster is defined as a collection of nodes and lines in the off-outage area that are supplied by only one available feeder. $\mathbb{X}(v,y^{(k)})$ determines the index of nodes that are in cluster $v$. $\mathbb{Y}(v,y^{(k)})$ denotes a set of configurations which provide the same optimal load pickup solutions for the nodes in cluster $v$. The formulation of $\mathbb{Y}$, optimality cuts and feasibility cuts are provided in section III.C.

Consider the simple network of Fig. 1 as an example. The switching operations shown in this figure are assumed to represent the optimal solution found in the master problem. Under the resulting configuration, the off-outage area is partitioned in two clusters. First cluster includes nodes 14, 15 and 16 that are restored from feeder-c through tie-switch T4. The second cluster includes nodes 11, 12 and 13 that is supplied from feeder-b through tie-switches T2 and T3.

As shown in Fig. 1., there is no path between two nodes belonging to the two different clusters except through the slack bus. We assume that the slack buses are effectively fixing the voltage set point at the top of each feeder (bus "Sub2" in Fig. 1.). Under this assumption, it can be said that the change of loading in one cluster (change of $x$ variables) does not change any state variable outside that cluster. Therefore, we solve a separate sub problem for each individual cluster. The aim is to break the computation burden of the inner level problem into several problems, which can be handled using different cores in parallel. The following MISOCP formulation models the sub problem for cluster $v$ at iteration $q$, given network configuration $y^{(q)}$.

$$\mathcal{L}_v(y^{(q)}) := \min_{x,z,u} \sum_{i \in \mathbb{X}(v,y^{(q)})} f_i(1-x_i) + \sum_{ij \in \mathbb{L}(v,y^{(q)})} d_{ij} z_{ij} \quad (8.a)$$

Subject to:
$$\mathbb{C}_1(x,y^{(q)}) \neq \emptyset \quad (8.b)$$
$$\mathbb{C}_2(x,y^{(q)},z,u) \neq \emptyset \quad (8.c)$$

where, $\mathbb{L}$ denotes the index of lines within cluster $v$.

It should be noted that the sub problem $\mathcal{L}_v^{(q)}$ incorporates only those variables that are related to the lines and nodes in cluster $v$. The objective function (8.a) is the minimization of the total energy not supplied in cluster $v$ as formulated in (1.a). Since the complicating variables $y$ are fixed, the big-M coefficients used in (4.a) do not appear in the sub problem which is relaxing again the computation burden of the inner level problem with respect to the original optimization problem $P$ given in (2).

$$UB^{(q)} = \sum_{v=1}^{v_{max}} \mathcal{L}_v(y^{(q)}) \tag{9}$$

According to (9), the summation of optimal objective values $\mathcal{L}_v(y^{(q)})$, associated with all the clusters, induce an upper bound $UB^{(q)}$ for the reliability optimal solution in the original optimization problem $P$ (2). The optimal solution of the sub problem $\mathcal{L}_v(y^{(q)})$ is also used to augment the feasibility cuts of the master problem as formulated in the first constraint of (5.d). In case that there is no feasible solution for the sub problem, a feasibility cut is generated and added to the master problem as formulated in the second constraint of (5.d).

*C. Generating the optimality and feasibility cuts*

In a given iteration k, if the sub problem $\mathcal{L}_v(y^{(k)})$ has a solution, say $x^{(k)}$, then clearly, $x = x^{(k)}$ is the optimal solution of P if y=$y^{(k)}$. We look for $\mathbb{Y}$ as a set of y-solutions, such that if we solve the sub problem while fixing y variable to any point $y'$ in this set, no better solution than $x^{(k)}$ can be found. It means that for any $y' \in \mathbb{Y}$, $\mathcal{L}_v(y') \geq \mathcal{L}_v(y^{(k)})$. This constraint is formulated in the first expression of (5.d) and denominated as an optimality cut.

If the sub problem $\mathcal{L}_v(y^{(k)})$ at iteration k is infeasible, then $\mathbb{Y}$ is defined as a set of $y$-solutions, such that if $y$ variable is fixed to any other point $y'$ in this set, the sub problem will be still infeasible. In other words, in order to break the infeasibility, variable $y$ should take values outside the set of $\mathbb{Y}$. This constraint is formulated in the second expression of (5.d) and referred as the feasibility cut.

Note that $\mathbb{Y}(v, y^{(k)})$ is associated with a given master problem solution $y^{(k)}$ and also with a cluster $v$. Also note that the solution of the master problem, say $y^{(k)}$, represents the network configuration and the solution of sub problem, say $x^{(k)}$, is the value of load pickup variables in cluster $v$. According to the definition of the optimality and feasibility cuts, in order to derive $\mathbb{Y}$, we should find network configurations $y'$ that lead to no better reliability solutions in cluster $v$, with respect to $\mathcal{L}_v(y^{(k)})$. The optimal solution of load pickup variables within cluster $v$ will not improve under configuration $y'$ with respect to the optimal values under configuration $y$ if the following conditions hold:

a) All the nodes in cluster $v$ that were connected to each other under configuration $y$ (identified by $\mathbb{X}(v, y^k)$), are still connected to each other under configuration $y'$.
b) The injection nodes that were supplying the nodes in cluster $v$ under configuration $y$ are the same as those under configuration $y'$,

Consider the test system of Fig. 1, as an example. As mentioned earlier, nodes 14, 15 and 16 are in the first cluster that is supplied by feeder-c through tie-switch T4 and by DG at node 5 through tie-switch T5. Tie-switches T4 and T5 are named as source lines. Source lines of cluster $v$ are defined as the lines at the border of cluster $v$ that are injecting active or-/and reactive power to the cluster. Considering the example shown in Fig. 1, assume that all the nodes in the first cluster should be restored except node 16, according to the solution of the sub problem $\mathcal{L}_v(y^{(k)})$. Now, by changing the configuration, it is obvious that the load at node 16 still cannot be restored if a) nodes 14, 15 and 16 are still connected to each other and if b) these nodes are supplied through the same source lines (tie-switches T4 and T5).

According to two conditions mentioned above, the set $\mathbb{Y}(v, y^{(k)})$ is expressed as in (10).

$$\mathbb{Y}(v, y^{(k)}) = \{y | \exists v' \leq v_{max}(y) :$$
$$\mathbb{X}(v, y^{(k)}) \subseteq \mathbb{X}(v', y), \Upsilon(v', y) \subseteq \Upsilon(v, y^{(k)})\} \tag{10}$$

where, $\Upsilon$ represent the index of source lines that are injecting power to cluster $v$.

In order to preserve the linearity of the optimality and feasibility cuts in terms of $y$ variables, the two constraints expressed in (10) are reformulated in (11) and (12).

$$\sum_{l \in \mathbb{L}(v, y^{(k)})} y_l = |\mathbb{X}(v, y^{(k)})| - 1 \tag{11}$$

$$y_l \leq y_l^{(k)} \quad : \quad \forall l \in \mathbb{L}(v, y^{(k)}) \tag{12}$$

The connectivity condition mentioned in condition a) is formulated in (11). This constraint enforces that the number of closed lines in a given cluster $v$ is equal to the total number of nodes in cluster $v$ minus one. This is the tree condition for cluster $v$. The tree condition ensures the network connectivity if it is radial [8]. This radiality condition is ensured for a given cluster $v$ using (5.b). Constraint (12) formulates condition b) that is mentioned above. This constraint ensures that the resource line $l$ of cluster $v$ that was open under configuration $y^{(k)}$ will stay open under any configuration $y \in \mathbb{Y}(v, y^{(k)})$.

According to the derived formulations for the set of $\mathbb{Y}(v, y^{(k)})$, the feasibility cut that was given in the second constraint of (5.d) is re-formulated in the following.

$$\sum_{l \in \mathbb{L}(v, y^{(k)})} y_l \leq |\mathbb{X}(v, y^{(k)})| - 2 \tag{13}$$

$$\sum_{l \in \mathbb{L}(v, y^{(k)})} y_l \geq 1 + \sum_{l \in \mathbb{L}(v, y^{(k)})} y_l^k \tag{14}$$

where, at least one of the conditions (13) or (14) must hold. This *either-or* constraint cannot be integrated into a convex model. Since in a convex optimization problem, all the constraints must hold. Therefore, this constraint should be further re-formulated according to the strategy given in Appendix II. Expressions (13) and (14) are the complements of (11) and (12), respectively. It should be mentioned that the complement of (12) means that at least one line $l \in \mathbb{L}(v, y^{(k)})$ exists such that $y_l \geq y_l^{(k)} + 1$. This constraint can be expresses as in (14).

Regarding the optimality cut, the first expression given in (5.d) can be translated into the following:

*If (11) and (12) are satisfied,*
*then $\sum_{i \in \mathbb{X}(v, y^{(k)})} f_i(1 - x_i) \geq \mathcal{L}_v(y^{(k)})$ must be satisfied.*

Therefore, the optimality cut is in the form of a conditional constraint. In order to be integrated in a convex optimization model, it is re-formulation as given in Appendix III.

*D. Modified Combinatorial Benders Algorithm*

The proposed decomposition approach for solving the distribution network restoration problem is described as follows:

1- Initialize iteration number ($q \leftarrow 1$), lower bound ($LB \leftarrow 0$), upper bound ($UB \leftarrow +\infty$), and set the convergence tolerance ($\varepsilon > 0$).
2- Solve the master problem to get the optimal network configuration $y^{(q)}$. Update the lower bound ($LB \leftarrow \max(LB, LB^{(q)})$), where $LB^{(q)}$ is given in (7).
3- Solve the sub problem for the obtained configuration $y^{(q)}$ and for each cluster $v$.
   a. If the optimization problem is feasible, find the optimal load pickup variables $x^{(q)}$ and augment the optimality cuts according to the first constraint of (5.d). Update the upper bound ($UB \leftarrow \min(UB, UB^{(q)})$), where $UB^{(q)}$ is given in (9).
   b. If the optimization problem leads to no feasible solution, augment the feasibility cuts according to the second constraint of (5.d).
4- Check for convergence :
   a. If $UB - LB \leq \varepsilon_{opt}$ or if the computation time is larger than $\varepsilon_{time}$, then terminate the algorithm and propose the best UB solution found so far as the solution of the problem.
   b. Else, update the iteration number ($q \leftarrow q + 1$), and return to step 2.

While the iterations of the proposed algorithm are evolving, the original solution space is gradually reduced by removing more combinations of binary variables. This is realized in a conservative way using the proposed optimality and feasibility cuts. Therefore, using the MCB approach, we might not be able to converge to the global optimal solution. However, as it will be illustrated in section IV.A, when the MCB algorithm converges, the best solution visited so far is close to the global optimal solution.

In order to end up with a tractable solution methodology in case of grids with realistic sizes, two stopping criteria are defined in the above-mentioned algorithm. According to this algorithm, we continue the running of iterations until the difference between the lower bound solution (LB) and the upper bound solution (UB) is lower than a threshold ($\varepsilon_{opt}$). In addition, we impose an additional threshold on the computation time ($\varepsilon_{time}$). In this regard, if the computation time is more than a threshold value, then the algorithm is stopped. The values of these thresholds ($\varepsilon_{opt}$ and $\varepsilon_{time}$) are determined based on the experience of DSO.

IV. NUMERICAL ANALYSIS AND DISCUSSION

In order to illustrate different features of the proposed solution algorithm for the restoration problem, two medium voltage networks are used shown in Fig. 2, and Fig. 3. In this paper, we study different scenarios. Scenarios I and II are applied on the test network of Fig. 2, whereas for scenarios III, the test network of Fig. 3 is used. For both test networks, the base power and energy values are assumed to be 1MW, and 1MWh, respectively. The minimum and maximum voltage magnitude limits are set, respectively, to 0.917 and 1.05 p.u [35]. The hourly profiles for different types of load patterns that are used in the both test networks are given in [36]. Two types of DGs are considered in this paper. Namely dispatchable and non-dispatchable. The dispatchable DGs, such as the diesel generators, are controlled to deliver the active and reactive power references that are set by DNO ahead of their operation. We consider also non-dispatchable DGs such as PV and wind generators, which are modeled as voltage-independent active power injection units with zero reactive power components. The forecast power injections of PV- and Wind-based DGs are derived from the real data reported in [37] and [38], respectively. In both test network, it is assumed that each node is equipped with a load breaker. All the line switches are assumed manually-controlled, whereas the load breakers are all assumed remotely-controlled. The time needed for the operation of each manually controlled and remotely-controlled switch are assumed 30 and 0.5 minutes, respectively. It is assumed that the critical loads that are shown with '⋈'. in Fig. 2 and Fig. 3, have the priority factors equal to 100 while the priority factors of other loads are equal to 1.

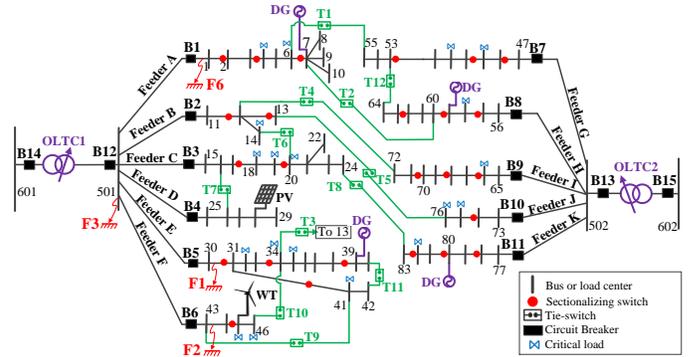

Fig. 2. The test network for test scenarios I and II [4].

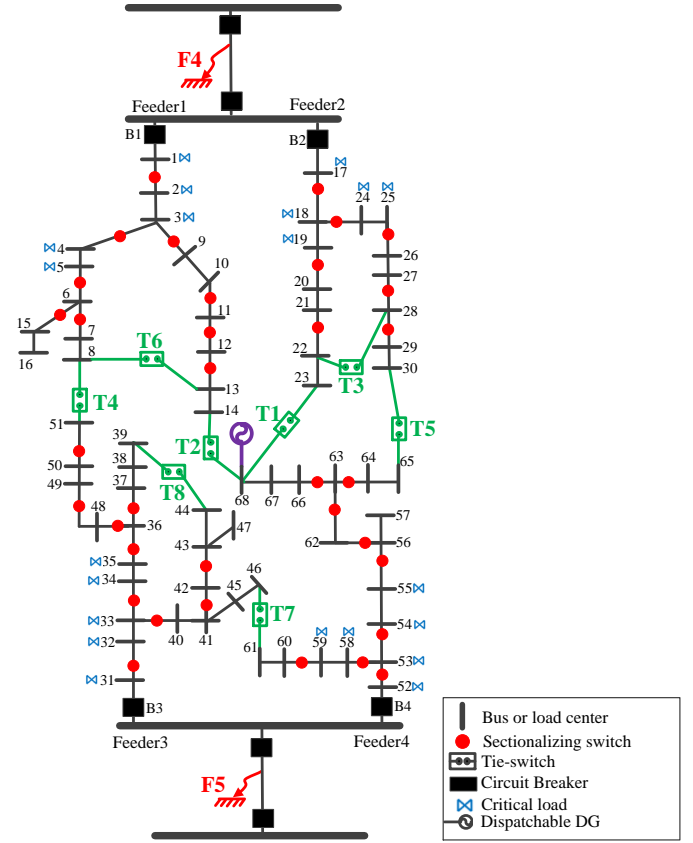

Fig. 3. The test network for test scenario III.

Table I. Comparison of restoration results obtained using IAO and MCB methods in scenarios I and II.

| Scenario | Solution Method | Reconfiguration Actions | Load Pickup Actions | $F^{re}$ (p.u.) | $F^{sw}$ (min) | Computation time (sec) |
|---|---|---|---|---|---|---|
| I (faults F1 and F2) | IAO | I. Open switch 38-39 and load breakers {33,34,37,41,42}, and close T11<br>II. Close T3 | - | 88 | 92.5 | 2.12 |
| I (faults F1 and F2) | MCB | I. Open switch 35-36 and load breakers {33,34,37,38,41,42}, and close T11<br>II. Close T3 | - | 100.1 | 93 | 9.58 |
| II (fault F3) | IAO | I. Open switches {1-2,33-34,31-41} and load breakers {2,3,4,5,6,12,13,14,35,36,37,38,39, 40,41}, and close {T3,T11}<br>II. Close{T1,T4} | III. Close load breaker 37 at t=12<br>IV. Close load breaker 13 at t=19 | 5.97e3 | 219 | 385.3 |
| II (fault F3) | MCB | I. Open switches {6-7,11-12,31-41} and load breakers {1,3,12,14,15,17,19,20,22,25,26,27, 28,29,30,31,32,34,35,36,38,43,44,46}, and close {T3,T7,T10}<br>II.Close{T1,T2,T5,T8} | III. Close load breakers {3,22,25,29} at t=19<br>IV. Close load breaker 13 at t=20 | 4.13e3 | 313.8 | 120 |

In order to show the functionality of the proposed solution approach, the restoration problem in case of each case study is solved using two approaches: I) the Integrated Analytical Optimization (IAO) method, and II) the Modified Combinatorial Benders (MCB) decomposition method proposed in section III. According to the IAO approach, the integrated optimization problem (1) is solved in one shot using an analytical solver. For this aim, the Branch-and-Bound method is used to relax the integrality constraints of the original optimization problem in an iterative way. In this regard, the best integer (valid) solution that is found at any step in the algorithm is called *incumbent* solution. The objective value of this incumbent is an upper bound for the optimal solution of the original minimization problem. At any step through the Branch-and-Bound search algorithm, there is also a lower bound, called the *best-bound* solution. This bound is obtained by taking the minimum of the optimal objective values of all the solutions obtained so far including the infeasible ones regarding the integrality constraints. The difference between the current upper and lower bounds is known as the *optimality gap*. It is said that the IAO approach converges to the proven optimality, when the optimality gap is less than the desired accuracy. This accuracy is tuned in this paper as 1e-10.

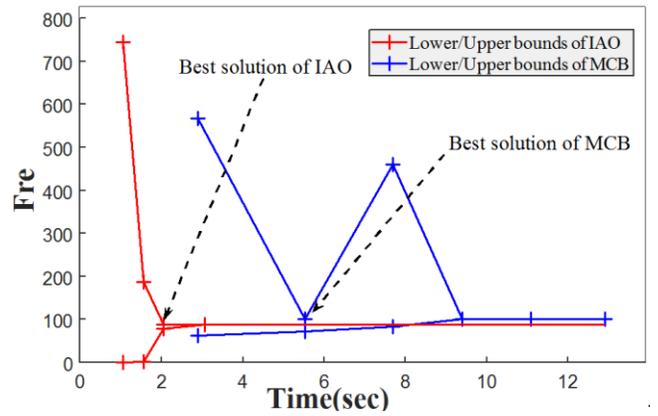

Fig. 4. The progress of the obtained solution using IAO and the MCB algorithms in scenario I.

The comparison of MCB and IAO methods is made, applying both methods on the same PC with an Intel(R) Xeon(R) CPU and 6 GB RAM, coded in Matlab/Yalmip environment and solved using Gurobi Optimizer 8.0. The restoration problem for all the test cases is considered as a multi-period optimization problem. The restorative period is assumed from 9:00 Am until 20:00 PM. The time step resolution is chosen to be 1 hour. We assume that the optimality threshold ($\varepsilon_{opt}$) and the computation time threshold ($\varepsilon_{time}$) are set to 0.01 p.u. and 2 minutes, respectively.

*A. Scenario I: a small-scale off-outage area*

The test network shown in Fig. 2 is a 11.4 kV balanced distribution network, which is based on a practical distribution system in Taiwan. It includes 2 substations, 11 feeders, 84 nodes, and 94 branches (incl. tie-branches). The detailed configuration data is given in [4]. Three dispatchable DGs on nodes 7, 39, and 80 have 2.8MW active and 3.0MVA apparent power capacities, while the DG on node 59 has 0.8MW active and 1.0MVA apparent power capacities. There exists also non-dispatchable DGs including a PV at node 28 and a Wind turbine at nodes 45.

In scenario I, the restoration problem is solved for the test network shown in Fig. 2, where two faults occur at the same time on lines 30-31 (fault F1) and 43-44 (fault F2). These two faults are isolated by opening line switches {B5, 30-31} and {B6, 44-45}, respectively. The resulting off-outage area includes feeder E except node {30} and feeder F except nodes {43 and 44}. The restoration solution obtained from IAO and MCB approaches are reported in Table I. In order to deploy this solution on the network, first, the "Reconfiguration Actions" must be implemented following the indicated order (I, II, etc.). Then, the "Load Pickup Actions" are deployed throughout the subsequent instants of the restorative period according to the schedule given in Table I. Along with these results, Table I provides the optimal values of different objective terms and the computation time.

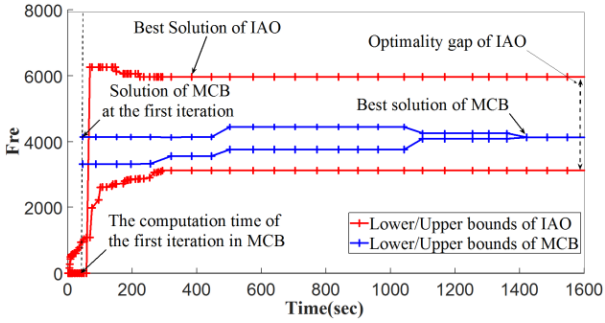

Fig. 5. The progress of the obtained solution using IAO and the MCB algorithms in scenario II.

The reliability objective term ($F^{re}$) is used to compare the quality of solutions provided by MCB and IAO methods[2]. In this regard, the quality margin of a restoration solution is defined as the difference between its reliability objective value and the global optimal value of the reliability objective term. As reported in Table I, the solution provided by the proposed MCB approach is 13.75% far from the global optimal solution, provided by IAO method.

The lower- and upper-bounds of the reliability objective term obtained using IAO and MCB approaches are plotted along their computation times in Fig. 4. In this regard, the lower and upper bounds of the MCB algorithm refer to the solutions provided, respectively, by the master and sub problem at each iteration. Whereas, for the IAO approach, each lower and upper bound correspond, respectively, to the best-bound and incumbent solutions found at a given iteration. It should be noted that in both methods only the upper bound solutions provide feasible solutions. As illustrated in Fig. 4, the best solution of the scenario I using the MCB algorithm is found at iteration 2 after 5.68 seconds.

*B. Scenario II: large-scale off-outage area in case of fault F3*

In this scenario, the restoration problem is studied in case of fault F3 at substation 501 in the test network of Fig. 2. The off-outage area includes the whole feeders A, B, C, D, E, and F. In this case, the optimization problem $P$ includes 529 binary variables including 23 reconfiguration variables ($y$) and 506 load pickup variables ($x$). As reported in Table I, it can be seen, the quality of the solution provided by the MCB method is 30.82% better than the one obtained with IAO method. In this scenario, IAO approach could not converge to the proven optimality. In this regard, the computation time that is given in Table I for IAO approach corresponds to the earliest time when IAO provides its best solution.

In scenario II, the computation time threshold is met and we have to stop after the iteration number 2. However, for the sake of illustration goals, we let the iterations to continue until LB and UB solutions converge. The results are shown in Fig. 5. This figure shows the same type of information as illustrated in Fig. 4 but for test scenario II. As it can be seen, the quality of the best-found solution of the MCB method after 2 minutes is only 0.05% far from the quality of the final solution found at the convergence. The functionality of the IAO approach for solving scenario II is also illustrated in Fig. 5. As it can be seen, the optimality gap could not be reduced below 18% in IAO approach.

Table II. Numerical results of test scenario III.

| Scenario | Solution Method | Reconfiguration Actions | Load Pickup Actions | $F^{re}$ (p.u.) | $F^{sw}$ (min) | Computation time (sec) |
|---|---|---|---|---|---|---|
| III.a (fault F4) | MCB | I. Open switches {6-7,11-12,25-26,27-28,28-29} and load breakers{1,3,5,19,24,25} II.Close{T1,T4} | III. Close load breakers {5,25} at 12:00 P.M. | 75.5 | 334 | 61.12 |
| III.b (fault F5) | MCB | I. Open switches {31-32,33-40,41-42,42-43,59-60,63-64} and load breakers {32,33,34,35,52,53,54,55,58} II.Close{T1,T4} | III. Close load breakers {35,58} at 12:00 P.M. | 114.6 | 245.5 | 65.76 |

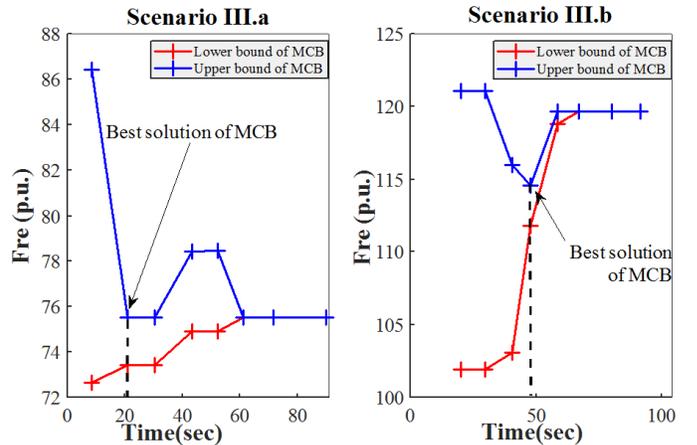

Fig. 6. The progress of the obtained solution using MCB algorithms in scenarios III.a and III.b.

*C. Scenario III: large-scale off-outage area in case of faults F4 and F5*

In order to further illustrate the performance of the proposed restoration algorithm, it is tested considering additional fault case studies. These are faults F4 and F5 represented in the test network of Fig. 3. This test network is based on a 11kV distribution network that is introduced in [39]. It includes 2 substations, 4 feeders, 70 nodes, and 76 branches (incl. tie-branches). One dispatchable DG with a capacity of 0.6 MW is installed on bus 68. The other type of DGs are PV-based DGs installed within LV networks at nodes 46, 47, and 61 with capacities of 0.6, 0.6 and 0.8 MW, respectively.

---

[2] Since $F^{re}$ has the largest weighting factor in the objective function, comparing the quality of the solution based on the overall objective value leads to the same result.

Table III. Checking network safety constraints for restoration solution obtained in different test cases

| Scenario | Fault | Min. voltage Margin (p.u.) | Min. current margin (A) |
|---|---|---|---|
| I | F1 and F2 | 0.0236 p.u. at node 31 at time 11:00 A.M | 9.716 A in line 84-11 at time 11:00 A.M |
| II | F3 | 0.0014 p.u. at node 29 at time 20:00 P.M. | 12.82 A in line 85-47 at time 18:00 P.M. |
| III.a | F4 | 0.0013 p.u. at node 25 at time 14:00 P.M. | 27.05 A in line 36-48 at time 13:00 P.M. |
| III.b | F5 | 0.0029 p.u. at node 31 at time 14:00 P.M. | 19.04 A in line 67-68 at time 14:00 P.M. |

In case of fault F4, feeders 1 and 2 will be in the off-outage area, whereas in case of fault F5, the off-outage area includes feeders 3 and 4. These two cases are presented, respectively, in scenarios III.a and III.b. The restoration problem contains 24 binary variables *y* and 186 binary variables *x* in case of scenario III.a; and 27 binary variables *y* and 234 binary variables *x* in case of scenario III.b.

The IAO approach fails to present even a single feasible restoration solution for scenarios III.a and III.b. The progress of the MCB algorithm in solving the restoration problem for these scenarios is depicted in Fig. 6. The numerical results corresponding to the best solution found until the convergence of the MCB algorithm are given in Table II. It can be seen that it takes only 21 and 48 seconds for the MCB to find these best solutions in case of scenarios III.a and III.b, respectively.

*D. Discussion*

The functionality of the proposed MCB with respect to the IAO method should be discussed separately for small scale-outage-areas such as in scenario I, and for large scale-outage-areas such as scenarios II and III. For small–scale problems, the IAO method provides the optimal solution within a short time. As mentioned in section III.C, the proposed MCB algorithm does not necessarily provide the global optimal solution (as for IAO). However, as shown in scenario I, the quality of its solution is not far from the global optimal solution.

The main advantage of the MCB method with respect to IAO method lies in large-scale optimization problems. Scenarios II and III illustrate this advantage. Actually, IAO approach could fail to converge to the optimal solution or even to provide a first feasible solution. As shown in scenario II, the best optimality gap obtained using IAO approach is significant and it means that the quality of the best feasible solution is poor. On the other hand, the MCB method found, within a very short time, a good-enough solution. There is no any unique and standard measure defining a good-enough restoration. However, it is obvious that the 0.05% quality margin that is obtained after the first iterations of the MCB algorithm in scenario II is acceptable and sufficient. In general, the IAO approach assigns the whole computational effort to finding the global optimal solution. If possible, this will be obtained for large-scale restoration problems after a long computation time. Whereas, the MCB method provides a good-enough solution at the first iterations thanks to the optimality cuts presented in section III.A. Then, through the subsequent iterations, it tries to improve the quality of the solution gradually. This characteristic is essential for the restoration problem, where an appropriate decision should be made in a very short time.

Table IV. Numerical restoration results obtained using the method of [24] in scenario II (fault F3 in the test network of Fig. 2).

| Solution Method | Reconfiguration Actions | Load Pickup Actions | $F^{re}$ (p.u.) | $F^{sw}$ (min) | Computation time (sec) |
|---|---|---|---|---|---|
| Method of [24]_try1 | I. Open switches {5-6,12-13,38-39,44-45} and load breakers {12,14,17,19,26,27,28,29,31,32,34,35,36,37,38,41,46}, and close {T3,T7,T10,T11} II.Close {T1,T2,T4,T5,T8} | III. Close load breaker 38 at 12:00 P.M. IV. Close load breakers {2,14,29} at 18:00 P.M. V. Close load breakers {35,36} at 19:00 P.M. | 3.32e3 | 1.11e3 | 28.04 |
| Method of [24]_try2 | I. Open switch 33-34 and load breakers {2,3,4,11,12,14,16,17,19,20,21,26,27,28,29,31,32,33,34,35,36,37,38,40,41,42,44,46}, and close {T3,T7,T9,T11} II.Close {T2,T4, T8} | III. Close load breakers {21,27,37} at 18:00 P.M. IV. Close load breakers {38,40} at 19:00 P.M. | 6.32e3 | 1.23e3 | 18.82 |

In order to check the network safety constraints for each scenario, the obtained restoration solution is deployed on the model of the corresponding test network implemented in Matlab environment. The voltage and current profiles along the time are derived using power flow simulations in Matlab/MATPOWER toolbox. Table III gives the representative numerical results out of these profiles for each scenario. These results include the minimum nodal voltage magnitude and minimum line current margins over, respectively, all the nodes and lines of the networks and over all the time steps during the restorative period. These results show that the network safety constraints are all respected and therefore confirm the feasibility of the obtained solutions. Moreover, it can be seen that according to each restoration solution, the network is operated very close to its capacity envelop (especially in terms of the minimum voltage limit). This illustrates that within the safe region of the network operation, the most possible loads are restored for each scenario.

*E. Comparison with other mathematical programming methods*

In this section, it is aimed to show the efficiency and superiority of the MCB with respect to the two mathematical programming methods proposed in [24] and [26]. In the first step, the MCB results in scenario II are compared with the results obtained from a mathematical formulation proposed in [24] for the restoration strategy.

In [24], the electrical safety constraints are integrated into the restoration problem using DistFlow formulation. In this regard, the multi-period restoration problem is formulated in [24] as a mixed-integer linear programming model. According to this formulation, the resulting topology of the network (i.e. line switching variables $y_{ij}$) could change at each time step.

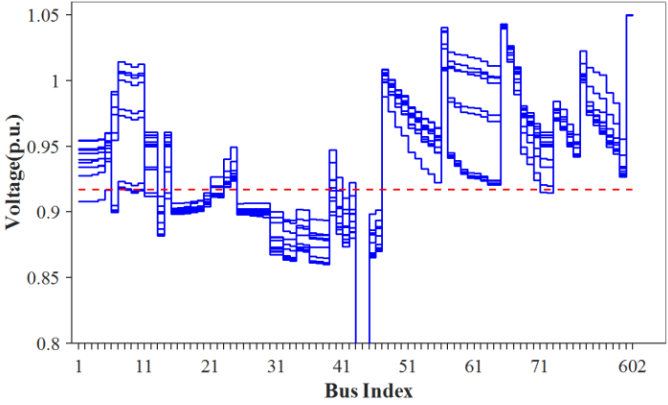

Fig. 7. The voltage magnitude profile at different times steps during the restorative period (blue lines) and lower voltage limit (red dotted line) according to the solution obtained from the method of [24]_try1 in scenario II (nodes 43 and 44 are left without any supply according to the results of Table IV).

In order to make a fair comparison, we force the line switching variables in [24] to not change with time, as suggested in the proposed MCB methodology. With this modification, the formulation of [24] is implemented in Yalmip/Matlab and solved using Gurobi for the test case of scenario II. The obtained numerical results are reported in the first row of Table IV (Method of [24]_try1).

The comparison of the $F^{re}$ value in Table IV with $F^{re}$ values reported in Table I for scenario II shows that the quality of the solution obtained using the method of [24] in try1 is better than the solution qualities of IAO and MCB approaches. However, since the DistFlow constraints are used in [24] instead of the AC power flow constraints, the obtained solution is infeasible regarding the minimum voltage limit. Fig. 7 confirms this infeasibility illustrating the results of a post power flow simulation for the obtained restoration solution. These results include the voltage magnitude profiles at different nodes and at different time steps during the restorative period. As it is shown, the lower voltage limit at some buses is violated at some time steps. The DistFlow approximation fails to guarantee the electrical safety constraints, since the network is operated close to its capacity envelop (i.e. current or/and voltage limits) during the restorative period.

In case where voltage and/or line power constraints are violated, it is suggested in [24] to impose conservative limits on the DistFlow constraints. In this regard, we replace the original lower voltage limit (0.917 p.u.) to 0.970 p.u. This is the smallest value for the lower voltage limit in the DistFlow constraints that can guarantee the feasibility of the obtained solution in scenario II.

Adding these conservative constraints, the restoration formulation of [24] is applied in try2 to solve the restoration problem in case of scenario II. The obtained results are shown in the second row of Table IV (Method of [24]_try2). This solution is feasible concerning all the electrical safety constraints. However, as it can be seen, the quality of the obtained solution (in terms of $F_{re}$) is 52.74% lower than the solution quality of the proposed MCB method reported in Table I.

Table V. Numerical restoration results in case of fault F6 in the test network of Fig. 2.

| Solution Method | | Switching Actions | $F^{re}$ (p.u.) | $F^{sw}$ (min) | Computation time (sec) |
|---|---|---|---|---|---|
| MCB Method | | I. Open switch 3-4  II. Close T2 | 4.4 | 60 | 4.04 |
| Method of [26] | Step1 | I. Open switch 6-7  II. Close {T1,T2} | 1.65e3 | 150 | 1.14 |
| | Step2 | Open load breakers {4,6} | | | |

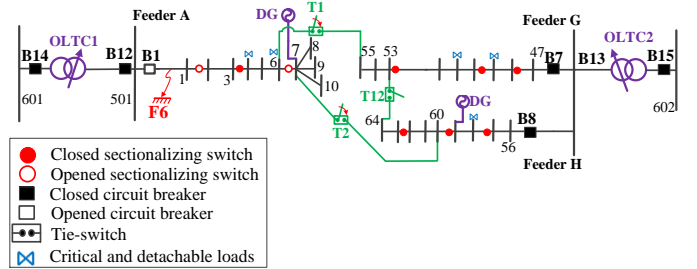

Fig. 8. The post-restoration configuration obtained in the first stage of the algorithm presented in [26].

This comparison clearly illustrates that it is not robust to simplify the AC power flow constraints in the restoration problem formulation with the corresponding DistFlow constraints. From the other hand, as the comparison of computation times related to the IAO and the method of [24] shows, the incorporation of AC power flow constraints increases the computation burden drastically. In order to make the restoration problem tractable in case of grids of realistic sizes while integrating the AC power flow constraints, the MCB method is proposed in this paper as a decomposition solution approach.

In the next step, the proposed MCB methodology is compared with the restoration methodology presented in [26]. First of all, it should be mentioned that the methodology of [26] is particularly suitable for unbalanced networks and in case of extreme fault cases, where there is no access to the upper grid for the network restoration (referred as islanded network restoration). In this section, the comparison is conducted considering only assumptions made in this paper. It means that we focus on balanced distribution networks and we do not consider the islanded restoration of the distribution network. For this comparison, it is assumed that a fault occurs at the top point of feeder A (fault F6) in the test network of Fig. 2. This fault is isolated by opening the feeder breaker B1 and the switch on line 1-2. All the parameters of this test case are similar to the ones of scenario I and scenario II except that all the nodes are not equipped with load breakers. In this case study, it is assumed that only critical loads that are shown with '⋈' can be detached from their nodes. It means that among all the nodes in the off-outage area, only nodes 4 and 6 are equipped with load breakers.

In order to make a fair comparison, we assume that load status variables ($\alpha_i$) will not change with time, as proposed by the authors of [26]. With this modification, we apply the developed MCB formulation on this test case. The numerical results are reported in Table IV. According to this solution, since the switch on line 3-4 is opened, the nodes 2 and 3 are left without any supply.

The restoration methodology of [26] is explained in section I. As mentioned, the authors of [26] propose to solve the restoration problem in two steps. In the first step, a heuristic approach is applied to find a suitable post-restoration topology for the network. This heuristic approach chooses a radial network configuration with the minimal diameter. The diameter of a tree is defined as the longest distance among all pairs of nodes in the network. The distance between a pair of nodes refers to the total impedance of lines on the shortest path between the two nodes. Applying this heuristic strategy to the test network of Fig. 2 in case of fault F6 results in the post-restoration configuration shown in Fig. 8.

According to the restoration algorithm presented in [26], in the second stage, the status of load breakers and the outputs of sources are determined while fixing the network topology to the one obtained in the first stage. The resulting optimization problem is solved in [26] using a relaxed semi-definite programming methodology. But in this paper, we solve this same optimization problem using Gurobi solver, which adopts the branch-and-bound method to handle integrality constraints. The optimal solution is to open all the load breakers in the off-outage area (i.e. at nodes 4 and 6). This leads to a reliability objective term equal to 1.65e3. The comparison of $F_{re}$ values given in Table IV shows that the quality of the solution obtained from the method of [26] is very far from the solution quality of the MCB method.

It should be noted that if there were no load breaker in the off-outage area (no detachable loads), there would be no feasible solution for the optimization problem in the second stage of the algorithm presented in [26]. These numerical results clearly highlights the limits of [26] mentioned in section I.

## V. Conclusion

The restoration is an NP-hard combinatorial optimization problem including three interdependent parts, namely, I) the reconfiguration problem, II) the load pickup problem, and III) the AC-Optimal Power Flow (AC-OPF) problem. This results in a huge and intractable problem especially considering a grid of realistic size in a multi-period problem. In order to relax the computation burden of the restoration problem, a two-stage decomposition approach is proposed in this paper named Modified Combinatorial Benders algorithm. In the outer level, the master problem solves a Mixed-Integer Linear optimization problem including the reconfiguration and the load pickup problems. The obtained configuration is fixed in the inner level and the load pickup variables are optimized subject to the AC-OPF constraints. The resulting sub problem is in the form of a Mixed-Integer Second-Order Cone Programming. This problem is broken down into several independent problems with smaller sizes. It makes the sub problem tractable in case of large-scale distribution networks. The solution of the sub problem is used to augment the feasibility or optimality cuts of the master problem. This algorithm is repeated through successive iterations until a solution with a desired level of optimality is obtained.

The superiority of the proposed decomposition approach with respect to the integrated approach is illustrated with different scenarios on two test distribution networks. The results indicate the functionality of the Modified Combinatorial Benders method in providing, within a short time, a good-enough solution for large-scale restoration problems. The future major work will expand the proposed decomposition method in order to account for the uncertainties in the forecast of load demands and DG production. These uncertainties will be incorporated to the optimization problem resulting in a multi-period stochastic restoration problem. As another step forward, we plan to extend the proposed formulation to the unbalanced distribution network according to the approach that we sketched in Appendix A.

## VI. Acknowledgment

The authors gratefully acknowledge the financial support of the Qatar Environment and Energy Research Institute (QEERI). The authors are also warmly thankful to prof. Jean-Yves Le Boudec (EPFL, Lausanne – Switzerland) for his review regarding the mathematical formulation.

## VIII. APPENDICES

### A. Appendix I: Discussion on the Extension of the Model to Unbalanced Grids

As mentioned in section I, the proposed approach is specifically derived for balanced distribution networks. A potential way to extend this strategy to generic unbalanced grids is described in this appendix. It is assumed that the distribution network is operated under unbalanced but still radial conditions during the restorative period. Moreover, we neglect the impact of non-transposed or partially transposed (asymmetrical) lines. Since the lines in MV distribution networks are relatively short and have small impedances, the impact of line asymmetry on the voltage unbalance may be ignored as it is negligible with respect to the impact of unbalanced load and generation [40]. In this regard, all the electrical state variables (not the switching variables) including voltage, current, and power flow variables are decomposed using well-known sequence transformation method. As a result, the unbalanced grid is decomposed into three symmetrical and balanced three-phase circuits. To each of these circuits, we apply the relaxed-OPF formulation as given in section II. Regarding the transformation of the voltage/current limits from phase domain to the sequence domain, we apply the methodology given in [16]. In this regard, we make a conservative assumption. We assume that the negative and zero sequence terms of voltage and current magnitudes are binding to their standard/normal limits. Therefore, the voltage and current limits associated with the positive sequence terms are derived *a priori*. The constraints regarding the OPF formulation of three sequences are integrated into the Master- and Sub-problems. Once the optimization problem is solved, we transform the obtained values of electrical state variables from sequence domain back into the phase domain.

B. *Appendix II: Convex formulation of either-or constraints*

As mentioned in section III.C, regarding the feasibility cut, at least one of two constraints (13) and (14) must hold. In order to integrate this either-or constraint into a convex optimization problem, binary variable $\psi$ is introduced subject to the following constraints:

$$\sum_{l\in \mathbb{L}(v,y^{(k)})} y_l \leq \left(\left|\mathbb{X}(v,y^{(k)})\right| - 2\right) + M_1 \psi \quad (15)$$

$$-M_2(1-\psi) + \sum_{l\in \mathbb{L}(v,y^{(k)})} y_l \geq 1 + \sum_{l\in \mathbb{L}(v,y^{(k)})} y_l^k \quad (16)$$

where, $M_1$ and $M_2$ are two positive and sufficiently-large numbers.

The auxiliary variable $\psi$ determines which of the two constraints (13) and (14) must hold. According to (15), if $\psi = 0$, (13) is imposed and (14) is relaxed. When $\psi = 1$, the situation is reversed. In both cases, one of the constraints is forced to be satisfied while the other constraint may also hold.

C. *Appendix III: Convex formulation of conditional constraints*

Let *"A implies B"* denotes a conditional constraint, where A and B are logical expressions. This is logically equivalent to state that *(A and ~ B)* is false, where $\sim B$ refers to the complement of B. The negation of *(A and ~ B)* is equivalent to *(~A or B)*.

In section III.C, a conditional expression was given for the optimality constraint, which is stated again in the following:

*If (11) and (12) are satisfied,*
*then $\sum_{i\in \mathbb{X}(v,y^{(k)})} f_i(1-x_i) \geq \mathcal{L}_v(y^{(k)})$ must be satisfied.*

According to the above explanation, this is logically equivalent to the following expression:

*(13) or (14) or $\sum_{i\in \mathbb{X}(v,y^{(k)})} f_i(1-x_i) \geq \mathcal{L}_v(y^{(k)})$ must hold*

As mentioned in section III.C, expressions (13) and (14) are the complements of (11) and (12), respectively. Therefore, the optimality cut can be translated into an either-or constraint and be reformulated according to the strategy given in section VIII.B.

D. *Appendix IV: Radiality Constraints*

In this appendix, the radiality constraints formulated in (1.b)-(1.d) are briefly explained. The details and illustrations of these formulations can be found in [32]. The orientation of line $ij$ with respect to the available tie-switch (as virtual sources) is determined by continuous variables $\beta_{ij}$ and $\beta_{ji}$. If the line is oriented from node $i$ to node $j$, variable $\beta_{ij}$ will be 1 and $\beta_{ji}$ will be zero and if the line is oriented from node $j$ to node $i$ variable $\beta_{ji}$ will be one and $\beta_{ij}$ will be zero.

In order to avoid possible isolated loops in cases where a DG exist in the off-outage area, (1.d) are added to the set of constraints. $M$ is a large and positive coefficient. To each line $ij$ with switch in the out-of-service area, two continuous flow variables $\Psi_{ij}$ and $\Psi_{ji}$ are assigned. They are associated with the binary variables $\beta_{ij}$ and $\beta_{ji}$, respectively. As formulated in the first two constraints of (1.d), for each line $ij$ with switch, at most one of the variables $\Psi_{ij}$ and $\Psi_{ji}$ gets a nonzero value depending on the line orientation that is identified with the variables $\beta_{ij}$ and $\beta_{ji}$. The third constraint of (1.d) formulates the flow (in a virtual commodity) balance equation for each node $i$, assuming that each node consumes a flow value equal to one. Finally, the last constraint of (1.d) implies that the total flows provided by all the available tie-switches must be equal to the total number of energized nodes. In [32] a mathematical proof is provided showing that the set of (1.b)-(1.d) constraints ensure a radial configuration without creating isolated loops supplied by DGs in islanded mode.